\documentclass[12pt]{article}
\usepackage[T2A]{fontenc}
\usepackage[utf8]{inputenc}
\usepackage[russian]{babel}
\usepackage{amssymb,bm,plain}
\begin{document}
\begin{plain}
\def\hyperrefs{\relax}
\font\eightrm=larm1000
\font\ninerm=larm1095
\font\nineit=lati1095
\font\ninebf=labx1095
\font\ninett=latt1095
\font\ninei=cmmi9 at 10.95pt
\font\ninesy=cmsy9 at 10.95pt
\font\nineex=cmex9 at 10.95pt
\font\tencmmib=cmmib10 at 12pt \skewchar\tencmmib='177
\font\sevencmmib=cmmib10 at 12pt \skewchar\sevencmmib='177
\font\fivecmmib=cmmib10 at 12pt \skewchar\fivecmmib='177
\newfam\cmmibfam
\textfont\cmmibfam=\tencmmib
\scriptfont\cmmibfam=\sevencmmib
\scriptscriptfont\cmmibfam=\fivecmmib
\skewchar\ninei='177
\skewchar\ninesy='60
\font\TITLE=labx1728
\def\fH{\mathfrak H}
\def\fD{\mathfrak D}
\def\fM{\mathfrak M}
\def\im{\mathop{\rm im}}

\def\No{\char 157}
\def\empty{}

\makeatletter
\def\Wo{{\mathpalette\Wo@{}}W}
\def\Wo@#1{\setbox0\hbox{$#1 W$}\dimen@\ht0\dimen@ii\wd0\raise0.65\dimen@%
\rlap{\kern0.35\dimen@ii$#1{}^\circ$}}

\catcode`\"=11
\def\k@vy{"}
\catcode`\"=13
\def"#1{\ifx#1<\char 190\relax\else\ifx#1>\char 191\relax\else%
	\ifx#1`\char 16\relax\else\ifx#1'\char 17\relax\else #1\fi\fi\fi\fi}

\def\hyp@rprov#1#2#3#4{\ifx\hyperrefs\undefined #1\else%
	\ifx#3\relax\else\fi #1%
	\ifx#3\relax\else\fi\fi}

\def\newl@bel#1#2{\expandafter\def\csname l@#1\endcsname{#2}}
\def\newbl@bel#1#2{\expandafter\def\csname bl@#1\endcsname{#2}}
\openin 11=\jobname .auxx
\ifeof 11
	\closein 11\relax
\else
	\closein 11
	\input \jobname .auxx
	\relax
\fi

\newcount\c@section
\newcount\c@subsection
\newcount\c@subsubsection
\newcount\c@equation
\newcount\c@bibl
\newcount\c@enum
\c@section=0
\c@subsection=0
\c@subsubsection=0
\c@equation=0
\c@bibl=0
\def\lab@l{}
\def\label#1{\immediate\write 11{\string\newl@bel{#1}{\lab@l}}\unskip}
\def\eqlabel#1{\rlap{$(\equation)$}\label{#1}}
\def\bibitem#1{\global\advance\c@bibl 1\par\noindent{\hyp@rprov{%
	\ninerm [\number\c@bibl]}{name=\k@vy cite:\number\c@bibl\k@vy}\relax{}~}%
	\immediate\write 11{\string\newbl@bel{#1}{\number\c@bibl}}}

\def\section#1{\global\advance\c@section 1
	{\par\vskip 3ex plus 0.5ex minus 0.1ex
	\rightskip=0pt plus 1fill\leftskip=0pt plus 1fill\noindent
	\hyp@rprov{{\bf\S\thinspace\number\c@section .~#1}}{%
	name=\k@vy sect:\number\c@section\k@vy}%
	\relax{}\par\penalty 10000\vskip 1ex plus 0.25ex}
	\gdef\lab@l{\number\c@section.:\number\pageno}
	\c@subsection=0
	\c@subsubsection=0
	\c@equation=0
}
\def\subsection{\global\advance\c@subsection 1
	\par\vskip 1ex plus 0.1ex minus 0.05ex\indent
	\hyp@rprov{{\bf\number\c@subsection.{}}}{%
	name=\k@vy sect:\number\c@section.\number\c@subsection\k@vy}\relax{}%
	\gdef\lab@l{\number\c@section.\number\c@subsection:0}
	\c@subsubsection=0
}
\def\subsubsection{\global\advance\c@subsubsection 1
	\par\vskip 1ex plus 0.1ex minus 0.05ex\indent%
	\hyp@rprov{{\bf\number\c@subsection.\number\c@subsubsection.{}}}{%
	name=\k@vy sect:\number\c@section.\number\c@subsection.%
	\number\c@subsubsection\k@vy}\relax{}%
	\gdef\lab@l{\number\c@section.\number\c@subsection.%
		\number\c@subsubsection:0}
}
\def\equation{\global\advance\c@equation 1
	\gdef\lab@l{\number\c@section.\number\c@equation:0}
	\hyp@rprov{{\rm\number\c@equation}}{name=\k@vy eq:\number\c@section(%
	\number\c@equation)\k@vy}\relax{}%
}
\catcode`\#=11
\def\sh@rp{#}
\catcode`\#=6
\def\ref@ref#1.#2:#3:{\def\REF@{#2}\ifx\REF@\empty%
	\hyp@rprov{{\rm \S\thinspace#1}}{href=\k@vy\sh@rp sect:#1\k@vy}0{rgb 0 0 1}%
	\else\hyp@rprov{\ifnum #1=\c@section {\rm #2}%
	\else {\rm \S\thinspace#1.#2}\fi}{href=\k@vy\sh@rp sect:#1.#2\k@vy}0{rgb 0 0 1}\fi
}
\def\ref@pageref#1:#2:{#2}
\def\ref@eqref#1.#2:#3:{\hyp@rprov{\ifnum #1=\c@section {\rm (#2)}\else%
	{\rm \S\thinspace#1$\,$(#2)}\fi}{href=\k@vy\sh@rp eq:#1(#2)\k@vy}0{rgb 0 0 1}%
}
\def\ref#1{\expandafter\ifx\csname l@#1\endcsname\relax
	{\bf ??}\message{^^J Reference #1 undefined!^^J}%
	\else\edef\mur@{\csname l@#1\endcsname :}%
	{\expandafter\ref@ref\mur@}\fi}
\def\pageref#1{\expandafter\ifx\csname l@#1\endcsname\relax
	{\bf ??}\message{^^J Reference #1 undefined!^^J}%
	\else\edef\mur@{\csname l@#1\endcsname :}%
	{\expandafter\ref@pageref\mur@}\fi}
\def\eqref#1{\expandafter\ifx\csname l@#1\endcsname\relax
	{(\bf ??)}\message{^^J Reference (#1) undefined!^^J}%
	\else\edef\mur@{\csname l@#1\endcsname :}%
	{\expandafter\ref@eqref\mur@}\fi}
\def\cite#1{\expandafter\ifx\csname bl@#1\endcsname\relax
	{\bf ??}\message{^^J Citation #1 undefined!^^J}%
	\else\hyp@rprov{{\bf\csname bl@#1\endcsname}}{href=\k@vy\sh@rp cite:%
	\expandafter\number\csname bl@#1\endcsname\k@vy}0{rgb 0 0 1}\fi}

\def\dfrac#1#2{{\displaystyle #1\over\displaystyle #2}}
\def\superalign#1{\tabskip=0pt\halign to\displaywidth{%
	\tabskip=0pt plus 1fil$\displaystyle ##$&%
	\tabskip=0pt\hss $\displaystyle ##{}$&%
	$\displaystyle {}##$\hss\tabskip=0pt plus 1fil&%
	\tabskip=0pt\hss ##\crcr #1\crcr}}
\long\def\enumerate#1{\c@enum=0\par\smallskip #1\par\smallskip}
\def\enumitem{\global\advance\c@enum 1\itemitem{$\number\c@enum^{\circ}$.\ }}
\makeatother

\def\proof{\par\medskip{\rm Д$\,$о$\,$к$\,$а$\,$з$\,$а$\,$т$\,$е$\,$%
	л$\,$ь$\,$с$\,$т$\,$в$\,$о.}\ }
\def\endproof{{}\hfill$\square$\par\smallskip}

\immediate\openout 11=\jobname.auxx

{\parindent=0cm УДК~517.984.42+519.218.7\par\vskip 2ex
\TITLE\rightskip=10pt plus 1fill\leftskip=10pt plus 1fill\noindent
Некоторые замечания об интегральных характеристиках винеровского процесса\par

\vskip 3mm\rm А.$\,$А.~Владимиров\footnote{}{\eightrm Работа поддержана
РФФИ, грант \No~10-01-00423.}\par}

\vskip 4ex
\section{Общие конструкции}\label{par:2}
\subsection\label{pt:1.1}
Обозначим через $\fH$ вещественное гильбертово пространство случайных величин
с конечным вторым моментом [{\bf \cite{GS:1977}}, Гл.~V, \S~1]. Рассмотрим
непрерывный в смысле нормы пространства $\fH$ центрированный гауссовский случайный
процесс $\xi:[0,1]\to\fH$, и обозначим через $K$ действующий в вещественном
пространстве $L_2[0,1]$ интегральный оператор, ядром которого выступает
ковариационная функция
$$
	{\mathcal K}(t,s)\rightleftharpoons\langle\xi(s),\xi(t)\rangle_{\fH}
$$
процесса $\xi$. Далее всегда будет предполагаться, что образ оператора $K$
является плотным в $L_2[0,1]$. Из этого предположения, в частности, вытекает
инъективность оператора $K$. Обозначим через $\fD$ пространство $\im K^{1/2}$,
снабжённое нормой
$$
	\|y\|_{\fD}\rightleftharpoons \|K^{-1/2}y\|_{L_2[0,1]},\leqno(\equation)
$$\label{eq:normD}%
а через $\fD^*$~--- двойственное к $\fD$ относительно $L_2[0,1]$
пространство, получаемое пополнением $L_2[0,1]$ по норме
$$
	\superalign{&\|y\|_{\fD^*}&\rightleftharpoons
		\sup\limits_{\|z\|_{\fD}=1}\langle y,z\rangle_{L_2[0,1]}&\cr
	&&=\langle Ky,y\rangle_{L_2[0,1]}^{1/2}.&[\eqref{eq:normD}]}
$$
Оператор $K$ при этом продолжается по непрерывности до изометрии $I:\fD^*\to\fD$
со свойством
$$
	(\forall y,z\in\fD^*)\qquad\langle y,Iz\rangle=\langle y,z\rangle_{\fD^*}.
	\leqno(\equation)
$$\label{eq:isomK}%

Заметим [{\bf \cite{GS:1977}}, Гл.~V, \S~1], что при любом выборе функции
$\varphi\in C[0,1]$ интеграл
$$
	\xi_{\varphi}\rightleftharpoons\int_0^1\varphi\cdot\xi\,dt
$$
корректно определён как бохнеровский (и даже римановский) интеграл
от \hbox{$\fH$-знач}\-ной функции. При этом выполняются не зависящие от выбора
$\varphi\in C[0,1]$ соотношения
$$
	\eqalign{\|\xi_{\varphi}\|^2_{\fH}&=\int_{[0,1]^2}{\mathcal K}(t,s)
		\varphi(s)\varphi(t)\,dtds\cr &=\|\varphi\|^2_{\fD^*},}
$$
означающие возможность продолжения соответствия $\varphi\mapsto\xi_{\varphi}$
по непрерывности до изометрического вложения пространства $\fD^*$ в пространство
$\fH$. Это наблюдение позволяет сопоставить каждому ортонормированному базису
$\{f_k\}_{k=0}^{\infty}$ пространства $\fD$ ортонормированное семейство
$\{\xi_k\}_{k=0}^{\infty}$ случайных величин $\xi_k\rightleftharpoons
\xi_{I^{-1}f_k}$, удовлетворяющее при любом выборе $\varphi\in\fD^*$
равенству
$$
	\superalign{(\equation)&\xi_{\varphi}&=\sum\limits_{k=0}^{\infty}
		\langle\varphi,f_k\rangle\,\xi_k.&[\eqref{eq:isomK}]}
$$\label{eq:preKL}%
В частности, независимо от выбора значения $t\in [0,1]$ дельта-функция
$\boldsymbol{\delta}_t$ с сосредоточенным в точке $t$ носителем заведомо
принадлежит пространству $\fD^*$ ввиду непрерывности ядра $\mathcal K$. Отвечающей
этой дельта-функции случайной величиной является значение процесса $\xi(t)$,
что автоматически влечёт [\eqref{eq:preKL}] справедливость разложения
$$
	\xi(t)=\sum\limits_{k=0}^{\infty}f_k(t)\,\xi_k,\leqno(\equation)
$$\label{eq:KL}%
известного как ряд Карунена--Лоева. При этом, ввиду компактности множества
дельта-функций в пространстве $\fD^*$, сходимость ряда \eqref{eq:KL} является
равномерной по параметру $t\in [0,1]$.

\subsection
Обозначим через $M$ действующий в пространстве $L_2[0,1]$ интегральный
оператор, ядром которого выступает ковариационная функция
$$
	{\mathcal M}(t,s)\rightleftharpoons\langle\xi^2(s),\xi^2(t)\rangle_{\fH}
$$
квадрата рассматриваемого процесса $\xi$. Обозначим через $\fM$ пополнение
пространства $L_2[0,1]$ по норме $\|y\|_{\fM}\rightleftharpoons\langle My,
y\rangle_{L_2[0,1]}^{1/2}$. Повторяя рассуждения из предыдущего пункта, легко
убеждаемся в возможности изометрического вложения пространства $\fM$
в пространство $\fH$ посредством соответствия
$$
	\tau_{\rho}\rightleftharpoons\int_0^1\rho\cdot\xi^2\,dt.\leqno(\equation)
$$\label{eq:101}%
Имеет место следующий факт.

\subsubsection\label{prop:2.tau}
{\it При любом выборе функции $\rho\in C[0,1]$ и ортонормированного базиса
$\{f_k\}_{k=0}^{\infty}$ пространства $\fD$ последовательность
$\{\tau_{\rho,n}\}_{n=0}^{\infty}$ случайных величин вида
$$
	\tau_{\rho,n}\rightleftharpoons\sum\limits_{k=0}^n\sum\limits_{l=0}^n
		r_{kl}(\rho)\cdot\xi_k\xi_l,\qquad r_{kl}(\rho)\rightleftharpoons
		\int_0^1\rho\cdot f_kf_l\,dt,
$$\label{eq:100}%
сходится в пространстве $\fH$ к интегралу \eqref{eq:101}.
}

\proof
Заметим, что для любой пары случайных величин $\zeta$, $\eta$ с центральным
нормальным совместным распределением выполняются соотношения
$$
	\eqalign{\|\zeta^2-\eta^2\|_{\fH}^2&=\langle (\zeta-\eta)^2,
		(\zeta+\eta)^2\rangle_{\fH}\cr
	&\leqslant \|(\zeta-\eta)^2\|_{\fH}\cdot\|(\zeta+\eta)^2\|_{\fH}\cr
	&=3\|\zeta-\eta\|_{\fH}^2\cdot\|\zeta+\eta\|_{\fH}^2.}
$$
Соответственно, последовательность квадратов частичных сумм ряда из правой части
равенства \eqref{eq:KL} равномерно по параметру $t\in [0,1]$ сходится
в пространстве $\fH$ к случайной величине $\xi^2(t)$. Последний факт немедленно
влечёт справедливость доказываемого утверждения.
\endproof

\subsection
Имеют место следующие два факта, связанные с понятиями {\it ядерного} оператора
и оператора {\it Гильберта--Шмидта} [{\bf \cite{GK:1965}}, Гл.~\hbox{III},
\S\S~8--9].

\subsubsection\label{prop:2.Sp}
{\it Пусть $R$~--- самосопряжённый ядерный оператор в арифметическом
гильбертовом пространстве $\ell_2$, а $\{\zeta_i\}_{i=0}^{\infty}$~---
система независимых стандартных нормальных случайных величин. Тогда двойной ряд
$$
	\dfrac{1}{\sqrt{3}+\sqrt{2}}\sum\limits_{\{i,j\}\in\mathbb N^2}
		R_{ij}\zeta_i\zeta_j\leqno(\equation)
$$\label{eq:sekv1}%
сходится в пространстве $\fH$, причём норма его суммы не превосходит ядерной
нормы оператора $R$.
}

\proof
Для любого конечного множества $\mathbb F\subset\mathbb N^2$ равенства $\|\zeta_i^2
\|_{\fH}=\sqrt{3}$ и ортонормированность системы $\{\zeta_i\zeta_j\}_{i>j}$
означают справедливость оценки
$$
	\eqalign{\left\|\sum\limits_{\{i,j\}\in\mathbb F} R_{ij}\zeta_i\zeta_j
		\right\|_{\fH}&\leqslant\sqrt{3}\cdot\sum\limits_{i=0}^{\infty}
		|R_{ii}|+\left(2\cdot\sum\limits_{\{i,j\}\in\mathbb N^2}
		|R_{ij}|^2\right)^{1/2}\cr
	&\leqslant(\sqrt{3}+\sqrt{2})\cdot\|R\|_1,}
$$
где через $\|R\|_1$ обозначена ядерная норма оператора $R$. Поскольку
последовательность $\{R_n\}_{n=0}^{\infty}$ операторов с матричными элементами
$$
	(R_n)_{ij}\rightleftharpoons\cases{R_{ij}&при $\max(i,j)\leqslant n$,\cr
		0&иначе}
$$
сходится по ядерной норме к оператору $R$, полученный результат означает
сходимость ряда \eqref{eq:sekv1} вместе с выполнением требуемой оценки нормы суммы.
\endproof

\subsubsection\label{prop:2.GS}
{\it Пусть $R$~--- самосопряжённый оператор Гильберта--Шмидта
в арифметическом гильбертовом пространстве $\ell_2$, удовлетворяющий тождеству
$R_{ii}\equiv 0$, а $\{\zeta_i\}_{i=0}^{\infty}$~--- система независимых
стандартных нормальных случайных величин. Тогда двойной ряд
$$
	\dfrac{1}{\sqrt{2}}\sum\limits_{\{i,j\}\in\mathbb N^2} R_{ij}\zeta_i\zeta_j
$$
сходится в пространстве $\fH$, причём норма его суммы совпадает с нормой
Гильберта--Шмидта оператора $R$.
}

\medskip
Справедливость данного утверждения немедленно вытекает из факта ортонормированности
системы случайных величин $\{\zeta_i\zeta_j\}_{i>j}$.


\section{Винеровский процесс}\label{par:3}
\subsection
Как хорошо известно, винеровский процесс определяется ковариационной функцией
$$
	{\mathcal K}(t,s)\equiv\inf(t,s),
$$
представляющей собой функцию Грина граничной задачи
$$
	\displaylines{-y''-\lambda y=0,\cr y(0)=y'(1)=0.}
$$
Соответственно, пространство $\fD$ в этом случае имеет вид
$$
	\fD=\{y\in W_2^1[0,1]: y(0)=0\},\qquad \|y\|_{\fD}=\|y'\|_{L_2[0,1]}.
$$
Ковариационная функция квадрата винеровского процесса имеет вид
$$
	{\mathcal M}(t,s)\equiv 2\inf(t,s)^2+ts.
$$
Непосредственным вычислением легко проверяется справедливость следующего
утверждения.

\subsubsection\label{prop:DM}
{\it Для любых функции $y\in L_2[0,1]$ и её первообразной $Y$ вида
$$
	Y(x)\equiv-\int_x^1 y\,dt\leqno(\equation)
$$\label{eq:prim}%
выполняются равенства
$$
	\leqalignno{\|y\|_{\fD^*}&=\|Y\|_{L_2[0,1]},&(\equation)\cr
	\label{eq:wiD}%
	\|y\|^2_{\fM}&=\left(\int_0^1Y\,dt\right)^2+4\int_0^1 t\,Y^2\,dt.
	&(\equation)}
$$\label{eq:wiM}
}

\noindent
Также имеют место следующие два факта.

\subsubsection
{\it Пространство $\fD^*$ непрерывным образом вложено в пространство $\fM$}
[\eqref{eq:wiD}, \eqref{eq:wiM}].

\subsubsection
{\it Любой элемент пространства $\fM$ представляет собой мультипликатор
из пространства $\fD$ в пространство $\fD^*$.
}

\proof
Зафиксируем произвольную функцию $\rho\in C[0,1]$ и отвечающую ей первообразную
$P$ вида
$$
	P(x)\equiv-\int_x^1\rho\,dt.\leqno(\equation)
$$\label{eq:P}%
Норма соответствующего самосопряжённого оператора умножения $R:\fD\to\fD^*$
допускает оценки
$$
	\superalign{&\|R\|&=\sup\limits_{\|y\|_{\fD}=1}
		\left|\int_0^1\rho\cdot y^2\,dt\right|&\cr
	&&=2\sup\limits_{\|y'\|_{L_2[0,1]}=1}\left|\int_0^1 P\cdot yy'\,
		dt\right|&\cr
	&&\leqslant 2\sup\limits_{\|y'\|_{L_2[0,1]}=1}\int_0^1
		\sqrt{t}|P|\cdot |y'|\,dt&\cr
	&&=2\left(\int_0^1 tP^2\,dt\right)^{1/2}&\cr
	&&\leqslant\|\rho\|_{\fM},&[\eqref{eq:wiM}]}
$$
что немедленно влечёт справедливость доказываемого утверждения.
\endproof

\subsection
Имеет место следующий факт.

\subsubsection\label{prop:Sp}
{\it Пусть обобщённая функция $\rho\in\fD^*$ определяет ядерный
мультипликатор типа $\fD\to\fD^*$. Пусть также последовательность
$\{\lambda_n\}_{n=0}^{\infty}$ перечисляет без повторений всевозможные
собственные значения граничной задачи
$$
	\displaylines{-y''-\lambda\rho y=0,\cr y(0)=y'(1)=0,}
$$
а последовательность $\{y_n\}_{n=0}^{\infty}$ составлена отвечающими
собственным значениям $\lambda_n$ нормированными в пространстве $\fD$
собственными функциями. Тогда выполняется равенство
$$
	\int_0^1\rho\cdot\xi^2\,dt=
		\sum\limits_{n=0}^{\infty}\lambda_n^{-1}\xi^2_{I^{-1}y_n}.
$$
}

\proof
Зафиксируем в пространстве $\fD$ ортонормированный базис $\{f_k\}_{k=0}^{%
\infty}$ вида
$$
	f_k(x)\equiv\dfrac{\sqrt{8}\,\sin [\pi\cdot(k+1/2)x]}{\pi\cdot(2k+1)}
$$
и свяжем с ним систему $\{\xi_k\}_{k=0}^{\infty}$ независимых стандартных
нормальных случайных величин $\xi_k\rightleftharpoons\xi_{I^{-1}f_k}$.
Сопоставим произвольной обобщённой функции $\rho\in\fD^*$ две случайные величины
$$
	\eqalign{\tau^{\bf d}_{\rho}&\rightleftharpoons
		\sum\limits_{i=0}^{\infty}r_{ii}(\rho)\cdot\xi_i^2,\cr
	\tau^{\bf n}_{\rho}&\rightleftharpoons\sum\limits_{\{i,j\}\in
		\mathbb N^2} r_{ij}(\rho)\cdot(1-\delta_{ij})\cdot\xi_i\xi_j,}
$$
где положено $r_{ij}(\rho)\rightleftharpoons\langle\rho,f_if_j\rangle$
и использован стандартный символ Кронекера $\delta_{ij}$. Соотношения
$$
	\superalign{&r_{kk}(\rho)&=\dfrac{8}{\pi^2\cdot(2k+1)^2}\cdot
		\langle\rho,\,\sin^2[\pi\cdot(k+1/2)x]\rangle\cr
		&&=-\dfrac{4}{\pi\cdot(2k+1)}\int_0^1 P\cdot
		\sin [\pi\cdot(2k+1)x]\,dx&[\eqref{eq:P}]}
$$
вкупе с равенством $\|P\|_{L_2[0,1]}=\|\rho\|_{\fD^*}$ гарантируют корректность
определения случайной величины $\tau^{\bf d}_{\rho}$ и выполнение оценки
$$
	\|\tau^{\bf d}_{\rho}\|_{\fH}\leqslant\sqrt{3}\cdot\|\rho\|_{\fD^*}.
$$
Кроме того, действующий из пространства $\fD$ в пространство $\fD^*$ оператор
умножения на функцию $\rho$ изометрично подобен действующему в пространстве
$L_2[0,1]$ интегральному оператору с ядром $-P(\sup(t,s))$, а потому
представляет собой оператор Гильберта--Шмидта. Это гарантирует [\ref{prop:2.GS}]
корректность определения случайной величины $\tau^{\bf n}_{\rho}$ и выполнение
оценок
$$
	\eqalign{\|\tau^{\bf n}_{\rho}\|_{\fH}&\leqslant
		\left(2\int_{[0,1]^2}P^2(\sup(t,s))\,dtds\right)^{1/2}\cr
	&\leqslant 2\cdot\|\rho\|_{\fD^*}.}
$$

Заметим теперь, что для любой функции $\rho\in C[0,1]$ выполняется
[\ref{prop:2.tau}] равенство
$$
	\tau_{\rho}=\tau^{\bf d}_{\rho}+\tau^{\bf n}_{\rho}.
	\leqno(\equation)
$$\label{eq:tauDN}%
Ввиду плотности линейного множества $C[0,1]$ в пространстве $\fD^*$ и непрерывности
характера зависимости каждой из случайных величин $\tau_{\rho}$,
$\tau^{\bf d}_{\rho}$ и $\tau^{\bf n}_{\rho}$ от параметра $\rho\in\fD^*$,
это означает выполнение равенства \eqref{eq:tauDN} и для всякой $\rho\in\fD^*$.
Иначе говоря, выполняется равенство
$$
	\tau_{\rho}=\sum\limits_{\{i,j\}\in\mathbb N^2}r_{ij}(\rho)\cdot
		\xi_i\xi_j.
$$
При этом в случае ядерности порождённого обобщённой функцией $\rho$ мультипликатора
последний может быть представлен в виде суммы сходящегося по ядерной норме
операторного ряда
$$
	\sum\limits_{n=0}^{\infty}\lambda_n^{-1}\langle I^{-1}y_n,\cdot\rangle\,
		I^{-1}y_n,
$$
что вкупе с равенствами
$$
	\superalign{&\sum\limits_{\{i,j\}\in\mathbb N^2}\langle I^{-1}y_n,
		f_i\rangle\cdot\langle I^{-1}y_n,f_j\rangle\cdot\xi_i\xi_j&
		=\xi_{I^{-1}y_n}^2&[\eqref{eq:preKL}]}
$$
как раз и означает [\ref{prop:2.Sp}] справедливость доказываемого утверждения.
\endproof

В случае неотрицательности мультипликатора $\rho$ утверждение \ref{prop:Sp} хорошо
известно и составляет~--- вместе с асимптотиками бесконечномерных \hbox{$\chi^2$-%
рас}\-пре\-де\-ле\-ний~--- основу теории малых уклонений винеровского процесса
(см., например, [\cite{Lif:1997}, \cite{Naz:2004}]). В общем случае оно указывает
возможные теоретико-вероятностные приложения спектральных асимптотик задачи
Штурма--Лиувилля с индефинитным весом (см., например, [\cite{VSh:2006:1},
\cite{VSh:2006:2}]).


\section{Пример неядерного мультипликатора}\label{par:A}
\subsection
Обозначим через $\rho_n$, где $n\geqslant 1$, обобщённые функции вида
$$
	\rho_n\rightleftharpoons\sum\limits_{k=1}^n\left(\boldsymbol{\delta}_{%
		(k-1/2)/n}-{\boldsymbol{\delta}}_{k/n}\right).
$$
Имеет место следующий факт.

\subsubsection\label{prop:A1}
{\it Пусть $\lambda_{m,n}$ есть \hbox{$m$-ое} снизу положительное собственное
значение граничной задачи
$$
	\displaylines{\rlap{$(\equation)$\label{eq:A1}}\hfill
		-y''-\lambda\rho_ny=0,\hfill\cr
		\rlap{$(\equation)$\label{eq:A2}}\hfill y(0)=y(1)=0.\hfill}
$$
Тогда независимо от выбора величины $\varepsilon>0$ при всех достаточно больших
$n\in\mathbb N$ выполняется неравенство $\lambda_{m,n}<2\pi m+\varepsilon$.
}

\proof
Введём обозначения
$$
	\omega_{k,n}(\lambda)\rightleftharpoons y_n(\lambda,k/n),\qquad
	\omega^{+}_{k,n}(\lambda)\rightleftharpoons y_n'(\lambda,[k/n]+0),
$$
где $y_n(\lambda,\cdot)$ есть решение уравнения \eqref{eq:A1} при начальных
условиях $y(0)=0$ и $y'(0)=1$. Как легко проверяется прямым вычислением,
для всех $k\in\{0,\ldots,n-1\}$ справедливы равенства
$$
	\pmatrix{\omega_{k+1,n}(\lambda)\cr \omega^{+}_{k+1,n}(\lambda)}=
		\left[1+\dfrac{1}{2n}\pmatrix{-\lambda&2-\lambda/(2n)\cr
		-\lambda^2&\lambda-\lambda^2/(2n)}\right]\cdot
		\pmatrix{\omega_{k,n}(\lambda)\cr\omega^{+}_{k,n}(\lambda)}.
$$
Это немедленно влечёт выполнение тождества
$$
	\pmatrix{\omega_{n,n}(\lambda)\cr\omega^{+}_{n,n}(\lambda)}=
		\left[1+\dfrac{1}{2n}\pmatrix{-\lambda&2-\lambda/(2n)\cr
		-\lambda^2&\lambda-\lambda^2/(2n)}\right]^n\cdot\pmatrix{0\cr 1},
$$
а тогда и равномерную внутри $\mathbb R$ сходимость функциональной
последовательности $\{y_n(\cdot,1)\}_{n=1}^{\infty}$ к функции $\omega\in
C(\mathbb R)$ вида
$$
	\omega(\lambda)\rightleftharpoons\cases{1&при $\lambda=0,$\cr
		\dfrac{2\sin(\lambda/2)}{\lambda}&иначе.}
$$
Последнее автоматически означает справедливость доказываемого утверждения.
\endproof

Отметим без доказательства, что несколько более детальное исследование позволило
бы установить, что \hbox{$m$-ые} снизу положительные собственные значения
граничных задач \eqref{eq:A1}, \eqref{eq:A2} при $n\to\infty$ в точности стремятся
к величине $2\pi m$.

\subsection
Обозначим через $\rho_{N,n}$, где $N,\,n\geqslant 1$, обобщённые функции вида
$$
	\rho_{N,n}\rightleftharpoons\sum\limits_{k=1}^n\left({\boldsymbol{\delta}}_{%
		(n+k-1/2)/(2^Nn)}-{\boldsymbol{\delta}}_{(n+k)/(2^Nn)}\right).
$$
Тривиальным образом [\ref{prop:DM}] имеют место равенства
$$
	\|\rho_{N,n}\|^2_{\fD^*}=2^{-N-1},
$$
означающие сходимость ряда
$$
	\rho_{\nu}\rightleftharpoons\sum\limits_{N=1}^{\infty}\rho_{N,\nu_N}
$$
в пространстве $\fD^*$ независимо от выбора последовательности $\{\nu_N\}_{%
N=1}^{\infty}$. Символы $\fD$ и $\fD^*$ здесь использованы для обозначения тех же
пространств, что и в \ref{par:3}.

Вариационные принципы для уравнения Штурма--Лиувилля с сингулярным весом (см.,
например, [\cite{VSh:2006:1}, \cite{VSh:2006:2}]) показывают, что \hbox{$m$-ое}
снизу положительное собственное значение граничной задачи
$$
	\displaylines{\rlap{$(\equation)$\label{eq:A3}}
		\hfill-y''-\lambda\rho_{\nu} y=0,\hfill\cr
		\rlap{$(\equation)$\label{eq:A4}}\hfill y(0)=y'(1)=0\hfill}
$$
заведомо мажорируется \hbox{$m$-ым} снизу положительным собственным значением
граничной задачи
$$
	\displaylines{-y''-\lambda\rho_{N,\nu_N} y=0,\cr
		 y(2^{-N})=y(2^{-N+1})=0.}
$$
Тем самым [\ref{prop:A1}], выбором последовательности $\{\nu_N\}_{N=1}^{\infty}$
можно добиться того, чтобы при всяком $m\geqslant 10$ соответствующее собственное
значение задачи \eqref{eq:A3}, \eqref{eq:A4} мажорировалось величиной $2\pi m
\ln m$. Это означает отсутствие безусловной сходимости ряда из величин, обратных
к собственным значениям задачи \eqref{eq:A3}, \eqref{eq:A4}, а потому и неядерность
определяемого функцией $\rho_{\nu}$ мультипликатора типа $\fD\to\fD^*$.

\penalty 10000\vskip 1ex
Автор выражает признательность А.$\,$И.~Назарову и И.$\,$А.~Шейпаку за ценные
обсуждения.

\vskip 0.75truecm
\centerline{\bf Список литературы}

\smallskip\ninerm
\bibitem{GS:1977} И.$\,$И.~Гихман, А.$\,$В.~Скороход. {\nineit Введение в теорию
случайных процессов}. М.: Наука, 1977.
\bibitem{GK:1965} И.$\,$Ц.~Гохберг, М.$\,$Г.~Крейн. {\nineit Введение в теорию
линейных несамосопряжённых операторов в гильбертовом пространстве}. М.: Наука,
1965.
\bibitem{Lif:1997} M.$\,$A.~Lifshits. {\nineit On the lower tail probabilities
of some random series}// Ann. Prob.~--- 1997.~--- V.~25, \No~1.~--- P.~424--442.
\bibitem{Naz:2004} А.$\,$И.~Назаров. {\nineit Логарифмическая асимптотика малых
уклонений для некоторых гауссовских процессов в \hbox{$L_2$-нор}\-ме
относительно самоподобной меры}// Записки науч. семин.~ПОМИ.~--- 2004.~---
Т.~311.~--- С.~190--213.
\bibitem{VSh:2006:1} А.$\,$А.~Владимиров, И.$\,$А.~Шейпак. {\nineit Самоподобные
функции в пространстве $L_2[0,1]$ и задача Штурма--Лиувилля с сингулярным
индефинитным весом}// Матем. сб.~--- 2006.~--- Т.~197, \No~11.~--- С.~13--30.
\bibitem{VSh:2006:2}А.$\,$А.~Владимиров, И.$\,$А.~Шейпак. {\nineit Индефинитная
задача Штурма--Лиувилля для некоторых классов самоподобных весов}// Труды
МИАН им.~В.$\,$А.~Стеклова.~--- 2006.~--- Т.~255.~--- С.~88--98.
\end{plain}
\end{document}